\documentclass[11pt]{article}

\usepackage{amsmath, amssymb, amscd}

\def\eqref#1{(\ref{#1})}
\newcommand{\goth}{\frak}

\newcommand{\arrow}{{\:\longrightarrow\:}}
\newcommand{\Z}{{\Bbb Z}}
\newcommand{\C}{{\Bbb C}}

\newcommand{\6}{\partial}
\def\1{\sqrt{-1}\:}

\newcommand{\calo}{{\cal O}}

\renewcommand{\tilde}{\widetilde}
\renewcommand{\bar}{\overline}
\renewcommand{\phi}{\varphi}
\renewcommand{\epsilon}{\varepsilon}
\renewcommand{\geq}{\geqslant}

\newcommand{\comment}[1]{{}}

\def\blacksquare{\hbox{\vrule width 4pt height 4pt depth 0pt}}
\def\endproof{\blacksquare}

\makeatletter

\@ifundefined{Bbb}
     {\newcommand{\Bbb}[1]{{\mathbb #1}}}%
{}

\newcommand{\ps@verbit}{%
  \renewcommand{\@oddhead}{%
          \scriptsize
          {Hypercomplex structures on K\"ahler manifolds}
          \hfil\tiny {M. Verbitsky, June 20, 2004}}
  \renewcommand{\@evenhead}{\@oddhead}
  \renewcommand{\@oddfoot}{\hfil\thepage\hfil}
  \renewcommand{\@evenfoot}{\@oddfoot}}
 
\newcounter{Mycounter}[section]
\newcounter{lemma}[section]
\renewcommand{\thelemma}{{Lemma \thesection.\arabic{lemma}}}
\newcommand{\lemma}{%
     \setcounter{lemma}{\value{Mycounter}}
     \refstepcounter{lemma}
     \stepcounter{Mycounter}
     {\bf \thelemma:\ }}

\newcounter{claim}[section]
\renewcommand{\theclaim}{{Claim \thesection.\arabic{claim}}}
\newcommand{\claim}{%
     \setcounter{claim}{\value{Mycounter}}
     \refstepcounter{claim}
     \stepcounter{Mycounter}
     {\bf \theclaim:\ }}

\newcounter{sublemma}[section]
\newcounter{corollary}[section]
\newcounter{theorem}[section]
\renewcommand{\thetheorem}{{Theorem \thesection.\arabic{theorem}}}
\newcommand{\theorem}{%
     \setcounter{theorem}{\value{Mycounter}}
     \refstepcounter{theorem}
     \stepcounter{Mycounter}
     {\bf \thetheorem:\ }}

\newcounter{conjecture}[section]
\newcounter{proposition}[section]
\renewcommand{\theproposition}
       {{Proposition \thesection.\arabic{proposition}}}
\newcommand{\proposition}{%
     \setcounter{proposition}{\value{Mycounter}}
     \refstepcounter{proposition}
     \stepcounter{Mycounter}
     {\bf \theproposition:\ }}

\newcounter{definition}[section]
\renewcommand{\thedefinition}
       {{Definition~\thesection.\arabic{definition}}}
\newcommand{\definition}{%
     \setcounter{definition}{\value{Mycounter}}
     \refstepcounter{definition}
     \stepcounter{Mycounter}
     {\bf \thedefinition:\ }}

\newcounter{example}[section]
\newcounter{remark}[section]
\renewcommand{\theremark}{{Remark \thesection.\arabic{remark}}}
\newcommand{\remark}{%
     \setcounter{remark}{\value{Mycounter}}
     \refstepcounter{remark}
     \stepcounter{Mycounter}
     {\bf \theremark:\ }}

\newcounter{problem}[section]
\newcounter{question}[section]
\renewcommand{\thequestion}{{Question \thesection.\arabic{question}}}
\newcommand{\question}{%
     \setcounter{question}{\value{Mycounter}}
     \refstepcounter{question}
     \stepcounter{Mycounter}
     {\bf \thequestion:\ }}

\@addtoreset{equation}{section}
\@addtoreset{footnote}{section}
\makeatother

\begin{document}

\begin{center}
{\LARGE\bf
Hypercomplex structures on \\[3mm] K\"ahler manifolds
}
\\[4mm]
Misha Verbitsky\footnote{ The author is
supported by EPSRC grant  GR/R77773/01 
and CRDF grant RM1-2354-MO02}
\\[4mm]

{\tt verbit@maths.gla.ac.uk, \ \  verbit@mccme.ru}
\end{center}

{\small 
\hspace{0.15\linewidth}
\begin{minipage}[t]{0.7\linewidth}
{\bf Abstract} \\
Let $(M,I)$ be a compact K\"ahler manifold
admitting a hypercomplex structure
$(M, I, J, K)$. We show that $(M, I, J, K)$
admits a natural HKT-metric. This is used to
construct a holomorphic symplectic 
form on $(M,I)$. 
\end{minipage}
}

{
\small
\tableofcontents
}

\section{Introduction}
\label{_Intro_Section_}


\subsection{Hypercomplex manifolds}

Let $(M, I, J, K)$ be a manifold equipped with an action of 
the quaternion
algebra ${\Bbb H}$ on $TM$. The manifold $M$ is called {\bf
hypercomplex} if the operators $I, J, K\in \Bbb H$ define
integrable complex structures on $M$. As Obata proved
(\cite{_Obata_}), this condition
is satisfied if and only if $M$ admits a torsion-free
connection $\nabla$ preserving the quaternionic action:
\[
\nabla I = \nabla J = \nabla K =0.
\]
Such a connection is called {\bf an Obata 
connection on $(M, I, J, K)$}. It is necessarily
unique (\cite{_Obata_}).

Hypercomplex manifolds were defined by C.P.Boyer
(\cite{_Boyer_}), who gave a classification of
compact hypercomplex manifolds for $\dim_{\Bbb H} M =1$.

If the Obata connection
$\nabla$, in addition, preserves a quaternionic Hermitian\footnote{
A metric $b$ is called {\bf quaternionic Hermitian} if 
\[ g(Ix, Iy) = g(Jx, Jy) = g(Kx, Ky) = g(x, y)\] for all $x, y \in TM$.}
metric $g$ on $M$, then $(M, I, J, K, g)$ is called {\bf hyperk\"ahler}.
This definition is equivalent to the standard one, see e.g. 
\cite{_Besse:Einst_Manifo_}.

\hfill

It is unknown precisely which complex manifold
admit hypercomplex structures. 

\hfill

\question
Consider a compact complex manifold $(M, I)$.
Describe the set of hypercomplex 
structures $(I, J, K)$ compatible with the 
given complex structure on $M$.

\hfill

A similar question about hyperk\"ahler structures is easily
answered by the Calabi-Yau theorem. Recall that a hyperk\"ahler
manifold is holomorphically symplectic. Indeed, consider
the 2-forms $\omega_J(\cdot, \cdot) = g(J\cdot, \cdot)$,
$\omega_K(\cdot, \cdot) = g(K\cdot, \cdot)$; then 
\begin{equation}\label{_Omega_Equation_}
\Omega := \omega_J + \1 \omega_K
\end{equation}
is a nowhere degenerate holomorphic
$(2,0)$-form on $(M, I)$ (\cite{_Besse:Einst_Manifo_}).
A converse result is implied by Calabi-Yau theorem:
a holomorphically symplectic compact 
K\"ahler manifold is necessarily 
hyperk\"ahler.

\hfill

\theorem\label{_Calabi_Yau_Theorem_}
Let $(M, I)$ be a compact holomorphically
symplectic manifold with a K\"ahler form $\omega$.
Then there exists a unique hyperk\"ahler metric
$g$ on $M$, with the same K\"ahler class
as $\omega$.

\hfill

{\bf Proof:} See \cite{_Besse:Einst_Manifo_}.
\endproof

\hfill

We have no similar description of complex manifolds admitting
a hypercomplex structure. In this paper we study the following
problem.

\hfill

\question
Let $(M, I)$ be a compact complex manifold 
of K\"ahler type\footnote{That is, admitting a 
K\"ahler metric.}. When $(M, I)$ admits a hypercomplex
structure?

\hfill

The following theorem gives an answer.

\hfill

\theorem\label{_main_intro_Theorem_}
Let $(M, I, J, K)$ be a compact hypercomplex
manifold. Assume that $(M, I)$ admits a K\"ahler
structure. Then $(M, I)$ is holomorphically symplectic.

\hfill

{\bf Proof:} In Subsection \ref{_HKT_cano_Subsection_} 
we deduce \ref{_main_intro_Theorem_}
from \ref{_canoni_trivi_Theorem_}, \ref{_HKT_exists_Theorem_}
and \ref{_HKT+Kah_=>Symple_Theorem_},
which are proven in Sections 
\ref{_CY_cano_tri_Section_}, \ref{_HKT_from_Kahler_Section_} and
\ref{_SS_on_HKT_Section_}. \endproof

\hfill

\remark
By Calabi-Yau theorem (\ref{_Calabi_Yau_Theorem_}), 
a holomorphically symplectic manifold admits a hyperk\"ahler
structure. However, the hypercomplex structure
$(M, I, J, K)$ on $M$ can {\it a priori} have a 
different nature. The manifold $(M, I, J, K)$
is hyperk\"ahler if and only if the Obata connection
$\nabla$ preserves a metric. However, if the holonomy
of $\nabla$ is non-unitarian, such a metric does not exist.

\hfill

\definition
Let $(M, I)$ be a compact holomorphically symplectic 
K\"ahler manifold, and $(M, I, J, K)$ a hypercomplex
structure on $(M,I)$. Then $(M, I, J, K)$
is called {\bf exotic} if $(M, I, J, K)$
is not hyperk\"ahler, that is, if the holonomy
of its Obata connection is not unitarian.

\hfill

We conjecture that exotic hypercomplex 
structures do not exist.

\subsection{HKT metrics and the canonical class}
\label{_HKT_cano_Subsection_}

Let $M$ be a hypercomplex manifold. A ``hyperk\"ahler with
torsion'' (HKT) metric on $M$ is a special
kind of a quaternionic Hermitian metric,
which became increasingly important in
mathematics and physics for the 
last 7 years.

HKT-metrics were
introduced by P.S.Howe and G.Papadopoulos (\cite{_Howe_Papado_})
and much discussed in physics literature since then.
For an excellent survey of these works written from  a mathematician's
point of view, the reader is referred to the paper of
G. Grantcharov and Y. S. Poon \cite{_Gra_Poon_}.

The term ``hyperk\"ahler metric with torsion'' is actually
misleading, because an HKT-metric is not hyperk\"ahler.
This is why we prefer to use the abbreviation ``HKT-manifold''.

Let $(M, I, J, K)$ be a hypercomplex manifold,
$g$ a quaternionic Hermitian form, and 
$\Omega$ the $(2,0)$-form on $(M, I)$ constructed
from $g$ as in \eqref{_Omega_Equation_}. The hyperk\"ahler
condition can be written down as $d\Omega=0$ 
(\cite{_Besse:Einst_Manifo_}).
The HKT condition is weaker:

\hfill

\definition
A quaternionic Hermitian metric is called an HKT-metric if
\begin{equation}\label{_HKT_intro_Equation_}
\6(\Omega)=0,
\end{equation}
where $\6:\; \Lambda^{2,0}_I(M) \arrow \Lambda^{3,0}_I(M)$
is the Dolbeault differential on $(M, I)$,
and $\Omega$ the $(2,0)$-form on $(M, I)$ constructed
from $g$ as in \eqref{_Omega_Equation_}. 

\hfill

It was shown in \cite{_Howe_Papado_}, 
\cite{_Gra_Poon_}, that this condition is
in fact independent from the choice of
the triple of complex structures 
$(I, J, K), IJ = -JI =K$ in ${\Bbb H}$.
In particular, we could replace the
hypercomplex structure 
$(M, I, J, K)$ with $(M, J, K, I)$.
We obtain the following trivial claim

\hfill

\claim\label{_HKT_remains_JKI_Claim_}
Let $(M, I, J, K)$ be a hypercomplex manifold,
$g$ a quaternionic Hermitian metric.
Consider $g$ as a quaternionic Hermitian metric
on a hypercomplex manifold $(M, J, K, I)$.
Then $g$ satisfies the HKT-condition
on $(M, I, J, K)$ if and only if
$g$ satisfies the HKT-condition
on $(M, J, K, I)$.

\endproof

\hfill

HKT-metrics play in hypercomplex geometry the same
role as the K\"ahler metrics play in complex geometry
(\cite{_Verbitsky:HKT_}). 

\hfill

The proof of \ref{_main_intro_Theorem_}
is split onto three steps, as follows. 

\hfill

\theorem\label{_canoni_trivi_Theorem_}
Let $(M, I, J, K)$ be a compact hypercomplex manifold.
Assume that $(M, I)$ admits a K\"ahler structure.
Then the ge exists a finite non-ramified covering
$\tilde M \arrow M$ such that the
canonical bundle of $(\tilde M, I)$ is trivial
as a holomorphic vector bundle.

\hfill

{\bf Proof:} See Section
\ref{_CY_cano_tri_Section_}. \endproof

\hfill

\theorem \label{_HKT_exists_Theorem_}
Let $(M, I, J, K)$ be a hypercomplex manifold. Assume that $(M, I)$
admits a K\"ahler metric $g$. Then $(M, I, J, K)$ admits an
HKT-metric $g_1$. Moreover, $g_1$ can be obtained 
by averaging $g$ with the $SU(2)$-action induced
by quaternions.

\hfill

{\bf Proof:} See Section \ref{_HKT_from_Kahler_Section_}. \endproof

\hfill

\theorem\label{_HKT+Kah_=>Symple_Theorem_}
Let $(M, I, J, K)$ be a compact hypercomplex manifold
admitting an HKT-metric. Assume that $(M, I)$
admits a K\"ahler structure. Assume, moreover, that there
exists a finite non-ramified cover $\tilde M \arrow M$ 
such that the canonical bundle $K(\tilde M, I)$ has a 
holomorphic trivialization.
Then $(M, I)$ is holomorphically symplectic.

\hfill

{\bf Proof:} See Section 
\ref{_SS_on_HKT_Section_}. \endproof

\hfill

\ref{_HKT+Kah_=>Symple_Theorem_} concludes the proof
of \ref{_main_intro_Theorem_}. Indeed, 
consider a compact hypercomplex manifold $(M, I, J, K)$,
and assume that $(M, I)$ admits a K\"ahler
metric. By \ref{_canoni_trivi_Theorem_}
the canonical class of $(M, I)$ is trivial, by
\ref{_HKT_exists_Theorem_} $(M, I, J, K)$ is HKT.
We arrive at assumptions of \ref{_HKT+Kah_=>Symple_Theorem_},
obtaining immediately that $(M, I)$ is holomorphically 
symplectic.

\hfill


\section{Calabi-Yau theorem and triviality of canonical bundle}
\label{_CY_cano_tri_Section_}


The following proposition is elementary.

\hfill

\proposition\label{_hc_c_1_zero_Proposition_}
Let $(M, I, J, K)$, $\dim_{\Bbb H}M =n$ be a hypercomplex manifold, and
\[ c_1(M, I) \in H^2(M, \Z)\] the first Chern class of $(M, I)$.
Then $c_1(M, I)=0$.

\hfill

{\bf Proof:} 
Let $SU(2)\subset {\Bbb H}^*$ be the group of 
unitary quaternions, acting on $TM$. A Riemannian
metric $g$ on $M$ is quaternionic Hermitian if and
only if $g$ is $SU(2)$-invariant. Taking an arbitrary
Riemannian metric and averaging over $SU(2)$, we 
obtain a quaternionic Hermitian metric. We proved the 
following trivial claim

\hfill

\claim
Let $M$ be a hypercomplex manifold. Then $M$ admits a quaternionic
Hermitian metric. 
\endproof

\hfill

Return to the proof of \ref{_hc_c_1_zero_Proposition_}.
To show that $c_1(M, I)=0$, we need to construct a 
continuous trivialization of the canonical bundle 
$K(M, I) = \Lambda^{2n,0}(M,I)$, where $2n = \dim_\C M$.
Let $g$ be a quaternionic Hermitian
metric on $M$, and 
\[ \Omega:= g(J \cdot, \cdot)+ \1 g(K \cdot, \cdot)\]
the corresponding non-degenerate $(2,0)$-form on $(M, I)$.
Then \[ \Omega^{n} \in \Lambda^{2n,0}_I(M)\]  is a non-degenerate
smooth section of the canonical bundle of
$\Lambda^{2n,0}_I(M)=K(M,I)$ of $(M,I)$.
Therefore, this bundle is topologically trivial.
This gives $c_1(M, I)=0$. \endproof

\hfill

The classification of K\"ahler manifolds with 
vanishing $c_1$ (\cite{_Bogomolov_}, 
\cite{_Beauville_}, \cite{_Besse:Einst_Manifo_})
easily implies the following result.

\hfill

\theorem \label{_trivi_cano_CY_Theorem_}
Let $(M, I)$ be a compact K\"ahler manifold with
\[ c_1(M, I)=0.\] Then there exists a finite non-ramified
covering $\tilde M \arrow M$ such that the 
canonical bundle $K(\tilde M,I)$ is trivial.

\endproof

\hfill

Combining \ref{_hc_c_1_zero_Proposition_}
and \ref{_trivi_cano_CY_Theorem_}, we obtain 
\ref{_canoni_trivi_Theorem_}.

\hfill

\remark
For a typical non-hyperkaehler
compact hypercomplex manifold $(M, I, J, K)$,
the complex manifold $(M, I)$ admits no K\"ahler metrics, and 
the Calabi-Yau theorem cannot be applied. The canonical bundle
$K(M,I)$ is trivial topologically by
\ref{_hc_c_1_zero_Proposition_}.
However, it is in most cases non-trivial 
as a holomorphic vector bundle, even if
one passes to a finite covering. It is 
possible to show that $K(M,I)$ is non-trivial
for all hypercomplex manifold $(M, I, J, K)$
such that $(M,I)$ is a principal toric fibration over
a base which has non-trivial canonical class; 
these include quasiregular Hopf manifolds and semisimple 
Lie groups with hypercomplex structure
constructed by D. Joyce (\cite{_Joyce_}).


\section{K\"ahler metrics and HKT metrics}
\label{_HKT_from_Kahler_Section_}


Let $(M, I, J, K)$ be hypercomplex manifold. Since $J$ and $I$
anticommute, $J$ maps $(p,q)$-forms on $(M,I)$
to $(q,p)$-forms:
\[ J:\; \Lambda^{p,q}_I(M) \arrow \Lambda^{q,p}_I(M).
\]

\definition
Let $\eta\in \Lambda^{2,0}_I(M)$ be a (2,0)-form on $(M,I)$.
Then $\eta$ is called {\bf $J$-real} if $J(\eta) = \bar\eta$,
and {\bf $J$-positive} if for any $x\in T^{1,0}(M,I)$,
$\eta(x, J(\bar x)) \geq 0$. We say that $\eta$ is 
{\bf strictly $J$-positive} if this inequality is strict
for all $x\neq 0$.

Denote the space of $J$-real, 
strictly $J$-positive  $(2,0)$-forms
by $\Lambda^{2,0}_{>0}(M,I)$.

\hfill

We need the following linear-algebraic lemma, which is
well known (see e.g. \cite{_Verbitsky_HKT_exa_}).

\hfill

\lemma\label{_2,0_forms_and_q_H_metrics_Lemma_}
Let $M$ be a hypercomplex manifold. Then 
$\Lambda^{2,0}_{>0}(M,I)$ is in one-to-one correspondence
with the set of quaternionic Hermitian metrics $g$ on $M$.
This correspondence is given by
\[ g \arrow  g(J \cdot, \cdot)+ \1 g(K \cdot, \cdot),\]
and the inverse correspondence by
\begin{equation}\label{_g_from_Omega_Equation_}
\Omega \arrow g(x,y):= \Omega(x, J(\bar y)).
\end{equation}

\endproof

\hfill

\lemma\label{_2,0_part_omega_J_Lemma_}
Let $(M,I,J,K)$ be a hypercomplex manifold, 
$g_1$ a Hermitian metric on $(M,J)$, $\omega_1= g_1(\cdot, J\cdot)$
the corresponding differential 2-form, and $\Omega_1$ the
$\Lambda^{2,0}_I(M)$-part of $\omega_1$.
Then $\Omega_1$ is strictly $J$-positive and $J$-real.

\hfill

{\bf Proof:} Since $\omega_1$ is a $(1,1)$-form on $(M,J)$,
we have $J(\omega_1)= \omega_1$. Therefore, $J(\Omega_1)=\bar\Omega_1$,
and $\Omega_1$ is $J$-real. 

Given $x\in T^{1,0}_I(M)$, $x \neq 0$,
the number $\Omega_1(x, J(\bar x))$ is real because 
$\Omega_1$ is $J$-real. On the other hand,
\[ \Omega_1(x, J(\bar x)) = \omega_1(x, J(\bar x))= g_1(x, \bar x)>0,\]
because $\Omega_1$ is a $(2,0)$-part of $\omega_1$ and
$x, J(\bar x)$ are $(1,0)$-vector fields.
We have shown that  $\Omega_1$ is strictly $J$-positive. This proves
\ref{_2,0_part_omega_J_Lemma_}. \endproof

\hfill

We also have the following trivial claim

\hfill

\claim\label{_6_Omega_1_=0_Claim_}
In assumptions of \ref{_2,0_part_omega_J_Lemma_},
let \[ \6:\; \Lambda^{p,q}_I(M) \arrow \Lambda^{p+1,q}_I(M)\]
denote the standard Dolbeault differential $\6$ on $(M,I)$.
Then $\6\Omega_1$ is the $(3,0)$-part of $d \omega_J$.
In particular, if $g_1$ is K\"ahler on $(M,J)$
then $\6\Omega_1=0$.

\hfill

{\bf Proof:} By definition, $\Omega_1$ is the $(2,0)$-part
of $\omega_J$, and $\6\Omega_1$ is the $(3,0)$-part of $d\Omega_1$.
\endproof

\hfill

\remark
Let $\phi$ be a K\"ahler potential for 
the K\"ahler form $\omega_1$ on $(M,J)$.
By Claim 2.3 of \cite{_Verbitsky_HKT_exa_},
on $(M,I)$ we have $\Omega_1= \6\6_J\phi$, where 
$\6_J = - J \circ \bar\6 \circ J$. 
The function $\phi$ satisfying $\Omega_1 = \6\6_J\phi$
for an HKT-form $\Omega_1$ is called 
{\bf an HKT-potential} for an HKT-form 
$\Omega_1$. 

\hfill

Now, let $(M,I,J,K)$ be a hypercomplex manifold, and
$g_1$ a K\"ahler metric on $(M,J)$. Consider the
form $\Omega_1\in \Lambda^{2,0}_I(M)$ constructed above.
Then $\Omega_1$ is strictly 
$J$-positive and $J$-real by \ref{_2,0_part_omega_J_Lemma_}, 
and hence corresponds
to a quaternionic Hermitian metric $g$ on $(M,I,J,K)$.
By \ref{_6_Omega_1_=0_Claim_}, $\6\Omega_1=0$, hence
$g$ is HKT. Doing all calculations explicitly, 
a reader can show that $g$ is obtained from $g_1$ 
by averaging over $SU(2)$ (we shall not use this 
claim). This proves \ref{_HKT_exists_Theorem_}.
Indeed, in assumptions of 
\ref{_HKT_exists_Theorem_} we are given a K\"ahler
metric on $(M,I)$, so the above argument
gives an HKT-metric on the hypercomplex manifold
$(M, J, K, I)$; this is equivalent to having an
HKT-metric on  $(M,I,J,K)$, as \ref{_HKT_remains_JKI_Claim_} 
implies. 


\section{Supersymmetry on HKT-manifolds with 
trivial canonical class}
\label{_SS_on_HKT_Section_}

Let $(M,I,J,K,g)$ be an HKT-manifold, and $K(M,I)$ 
its canonical class. Using the quaternionic Hermitian metric $g$
we trivialize the canonical class 
by a smooth non-degenerate section
as in \ref{_hc_c_1_zero_Proposition_}.
Let $K^{1/2}$ be the square root
of the canonical bundle corresponding to this
trivialization. Writing $K(M,I)$ as a trivial
bundle with the Chern connection
$\nabla_{triv}+ \theta$, we define 
$K^{1/2}$ as a trivial bundle
with the conection $\nabla_{triv}+ \frac 1 2\theta$.
This connection is clearly induced by a holomorphic structure
on $K^{1/2}$, and $K^{1/2}\otimes K^{1/2}$ is isomorphic
to $K$ as a holomorphic line bundle.

\hfill

In \cite{_Verbitsky:HKT_} we proved the following theorem, which
is implied by an analogue of the 
Lefschetz\--type $\goth{sl}(2)$-action in the HKT setting.

\hfill

\theorem\label{_SS_on_HKT_Theorem_}
Let $(M,I,J,K)$ be a compact HKT-manifold, $\dim_{\Bbb H} M=n$,
and $K^{1/2}$ the square root of a canonical
bundle $K(M,I)$ constructed as above. Consider the
Dolbeault class $[\bar\Omega]\in H^{0,2}_{\bar\6}(M,I)=H^2(\calo_{(M,I)})$
of $\bar\Omega$, where $\Omega\in \Lambda^{2,}_I(M)$ is 
the HKT-form of $M$, and let 
\begin{equation}\label{_Lefschetz_on_cohomo_HKT_Equation_}
H^{l}(K^{1/2}) \stackrel{\bigwedge [\bar\Omega]^{n-l}}\longrightarrow H^{2n-l}(K^{1/2})
\end{equation}
be the corresponding multiplicative map on the 
holomorphic cohomology of $K^{1/2}$. Then
\eqref{_Lefschetz_on_cohomo_HKT_Equation_} is an isomorphism.

\hfill

{\bf Proof:} In \cite{_Verbitsky:HKT_} it was shown that
the natural operator \[ L_{\Omega}:\; H^{l}(K^{1/2})\arrow H^{l+2}(K^{1/2})\]
belongs to an $\goth{sl}(2)$-triple. This is used in 
\cite{_Verbitsky:HKT_} to 
obtain \ref{_SS_on_HKT_Theorem_} in the same way as
one obtains a similar result for the 
cohomology of a K\"ahler manifold. \endproof

\hfill

When $K(M,I)$ is a trivial holomorphic
bundle, $K^{1/2}$ is also a trivial bundle. We obtain that
\begin{equation}\label{_Omega^n_minipage_Equation_}
\begin{minipage}[m]{0.8\linewidth}
when $K(M,I)$ is trivial, 
$[\Omega]^n$ is a generator 
of $H^{2n}(K^{1/2})\cong H^0(K^{1/2})^*=\C$
\end{minipage}
\end{equation}
(the last isomorphism is provided by the Serre's duality,
using the triviality of the canonical bundle). Now we can prove
\ref{_HKT+Kah_=>Symple_Theorem_}.

Let $(M,I,J,K)$ be a compact HKT-manifold, with
$\tilde M$ a non-ramified finite covering of $M$
with the canonical bundle $K(\tilde M,I)$ trivial.
Assume that $(M,I)$ admits a K\"ahler metric. 
By Calabi-Yau theorem (\cite{_Yau:Calabi-Yau_}),
$(M,I)$ admits a Ricci-flat K\"ahler metric $h$. Let
$\bar \Omega_h\in \Lambda^{0,2}_I(M)$
be a harmonic representative of the cohomology
class $[\bar\Omega]\in H^{0,2}_{\bar\6}(M,I)$ under
$h$. Since $(M,I)$ is K\"ahler, the harmonic
 $(2,0)$-form $\Omega_h$ is holomorphic.
By Bochner-Lichnerowicz theorem
(\cite{_Besse:Einst_Manifo_}),
this implies 
\begin{equation}\label{_Omega_h_parallel_Equation_}
\nabla_h\Omega_h=0,
\end{equation}
where $\nabla_h$ is the Levi-Civita connection
of $h$ (this is true for any
holomorphic form $\Omega_h$ on a Ricci-flat 
compact K\"ahler manifold). Let $\tilde \Omega_h$, $\tilde Omega$
be $\Omega_h$, $\Omega$ lifted to $\tilde M$.  
By \eqref{_Omega^n_minipage_Equation_},
$\tilde \Omega^n$, and hence
$\tilde \Omega^n_h$, represents non-zero 
class in cohomology of $(\tilde M,I)$. 
This implies $\Omega^n_h\neq 0$. By
\eqref{_Omega_h_parallel_Equation_},
we also have $\nabla_h\Omega_h^n=0$,
hence $\Omega^n_h$ trivializes 
$\Lambda^{2n,0}_I(M)$. We obtain that
$\Omega_h$ is a non-degenerate 
holomorphic symplectic form on $(M,I)$.
This proves \ref{_SS_on_HKT_Theorem_}.
We finished the proof of \ref{_main_intro_Theorem_}.

\hfill

{\bf Acknowledgements:}
This paper appeared as a result of a very rewarding
colloboration with S. Alesker. I am also grateful
to D. Kaledin and D. Kazhdan for interesting discussions,
and to S. Alesker for pointing out errors in the early
version of this manuscript. 

\hfill

{\scriptsize

\hfill

\noindent {\sc Misha Verbitsky\\
University of Glasgow, Department of Mathematics, \\
15 University Gardens, Glasgow G12 8QW, Scotland}, \\
{\sc  Institute of Theoretical and
Experimental Physics \\
B. Cheremushkinskaya, 25, Moscow, 117259, Russia }\\
\tt verbit@maths.gla.ac.uk, \ \  verbit@mccme.ru 
}

\end{document}